\documentclass[12pt]{article}       % DEFAULT font size is 10pt.
\usepackage{amsfonts,amsmath,amssymb}
\newtheorem{Def}{Definition}
\newtheorem{Th}{Theorem}

\newtheorem{lem}[Th]{Lemma}

\newtheorem{Conjecture}{Conjecture}

\newenvironment{Pf}
{\noindent\textbf{Proof.}\ \ }{\hfill $\Box$ \medskip}

\title{Equitable colorings of complete multipartite graphs}

\author {Keaitsuda Maneeruk Nakprasit\\ {\small\em Department of Mathematics, Faculty of Science, Khon Kaen University, 40002, Thailand }\\
{\small\em E-mail address: kmaneeruk@hotmail.com}
\and Kittikorn Nakprasit \footnote{Corresponding Author} \\
{\small\em Department of Mathematics, Faculty of Science, Khon Kaen University, 40002, Thailand }\\
{\small\em E-mail address: kitnak@hotmail.com}}
\date{}

\begin{document}

\maketitle
\begin{abstract} 
A $q$-\emph{equitable coloring} of a graph $G$ is a proper $q$-coloring  
such that the sizes of any two color classes differ by at most one. 
In contrast with ordinary coloring, a graph may have an equitable
$q$-coloring  but has no equitable $(q+1)$-coloring. 
The \emph{equitable chromatic threshold} 
is the minimum $p$ such that $G$ has an equitable
$q$-coloring for every $q\geq p.$ 

In this paper, we establish the notion of $p(q: n_1,\ldots, n_k)$ 
which can be computed in linear-time and prove the following. 
Assume that $K_{n_1,\ldots,n_k}$ has an  equitable $q$-coloring.   
Then $p(q: n_1,\ldots, n_k)$ is the minimum $p$ such that 
$K_{n_1,\ldots,n_k}$  has an equitable $r$-coloring for each $r$ 
satisfying $p \leq r \leq q.$ 
Since $K_{n_1,\ldots,n_k}$ has 
an equitable $(n_1+\cdots+n_k)$-coloring, 
the equitable chromatic threshold of $K_{n_1,\ldots,n_k}$ 
is $p(n_1+\cdots+n_k: n_1,\ldots, n_k).$ 

We find out later that the aforementioned immediate consequence 
is exactly the same as the formula of Yan and Wang \cite{YW12}. 
Nonetheless, the notion of $p(q: n_1,\ldots, n_k)$ 
can be used for each $q$ 
in which $K_{n_1,\ldots,n_k}$ has an equitable $q$-coloring 
and the proof presented here is much shorter. 
\end{abstract}
%%%%%%%%%%%%%%%%%%%%%%%%%%%%%%%%%%
% \begin{keyword}
% Equitable coloring \sep $k$-tree-coloring \sep Vertex k-arboricity 
% \sep Complete multipartite graph
% \end{keyword}

\section{Introduction}
Throughout this paper, all graphs are finite, undirected, and 
simple. We use  $V(G)$ and $E(G),$ respectively, 
to denote the vertex set and edge set of a graph $G.$ 
Let $K_{n_1,\ldots,n_k}$ be a complete $k$-partite graph in which  
partite set $X_i$ has size $n_i$ for $1 \leq i \leq k.$ 
Let $K_{k*n}$ denote a complete $k$-partite set with each partite set has size $n.$ 

An \emph{equitable $k$-coloring} of a graph is a proper vertex
$k$-coloring such that the sizes of every two color classes differ
by at most $1.$ 

It is known \cite{GareyJohnson} that determining if a planar graph
with maximum degree $4$ is $3$-colorable is NP-complete. For a given
$n$-vertex planar graph $G$ with maximum degree $4,$ let $G'$ be the
graph obtained from $G$ by adding $2n$ isolated vertices. Then $G$
has $3$-coloring if and only if $G'$ has an equitable $3$-coloring.
Thus, finding the minimum number of colors needed to color a graph
equitably even for a planar graph is an NP-complete problem.

Hajnal and Szemer\'edi~\cite{HS} settled a conjecture of Erd\H os by
proving that  every graph $G$ with maximum degree at most $\Delta$
has an equitable $k$-coloring for every $k\geq 1+\Delta.$ 
This result is now known as Hajnal and Szemer\'edi Theorem. 
Later, Kierstead and
Kostochka~\cite{KK08} gave a simpler proof of Hajnal and Szemer\'edi
Theorem. The bound of the
Hajnal-Sz{e}mer\' edi theorem is sharp, but it can be improved for
some important classes of graphs. In fact, Chen, Lih, and
Wu~\cite{CLW94} put forth the following conjecture.

\begin{Conjecture} \label{ConjLW}
Every connected graph $G$ with maximum degree $\Delta\geq 2$ has an
equitable coloring with $\Delta$ colors, except when $G$ is a
complete graph or an odd cycle or $\Delta$ is odd and
$G=K_{\Delta,\Delta}.$
\end{Conjecture}

Lih and Wu~\cite{LW} proved the conjecture for bipartite graphs.
Meyer \cite{M} proved that every forest with maximum degree $\Delta$
has an equitable $k$-coloring for each $k \geq 1+\lceil
\Delta/2\rceil $ colors. This result implies the conjecture holds
for forests.  Yap and Zhang~\cite{YZ1} proved that the
conjecture holds for outerplanar graphs. Later Kostochka~\cite{Ko}
improved the result  by proving that every
outerplanar graph with maximum degree $\Delta$ has an equitable
$k$-coloring for each $k \geq 1+\lceil \Delta/2\rceil.$ 

In~\cite{ZY98}, Zhang and Yap essentially proved the conjecture
holds for planar graphs with maximum degree at least $13.$ Later
Nakprasit~\cite{Nak12} extended the result to all planar graphs with
maximum degree at least  $9.$ 
Some related results are about planar graphs without some restricted cycles 
~\cite{LiBu09, NN12, ZhuBu08}. 

Moreover, the conjecture has been confirmed for other classes of graphs, 
such as graphs with degree at most 3~\cite{CLW94, CY12}  
and series-parallel graphs \cite{ZW11}. 

In contrast with ordinary coloring, a graph may have an equitable
$k$-coloring  but has no equitable $(k+1)$-coloring. 
For example, $K_{7,7}$ has an equitable $k$-coloring
for $k=2,4,6$ and $k \ge 8$, but has no equitable $k$-coloring for
$k=3,5$ and $7$. This leads to the definition of the 
\emph{equitable chromatic threshold} which is 
the minimum $p$ such that $G$ has an equitable
$q$-coloring for every $q\geq p.$ 

In this paper, we establish the notion of $p(q: n_1,\ldots, n_k)$ 
which can be computed in linear-time and prove the following. 
Assume that $K_{n_1,\ldots,n_k}$ has an  equitable $q$-coloring.   
Then $p(q: n_1,\ldots, n_k)$ is the minimum $p$ such that 
$K_{n_1,\ldots,n_k}$  has an equitable $r$-coloring for each $r$ 
satisfying $p \leq r \leq q.$ 
Since $K_{n_1,\ldots,n_k}$ has 
an equitable $(n_1+\cdots+n_k)$-coloring, 
the equitable chromatic threshold of $K_{n_1,\ldots,n_k}$ 
is $p(n_1+\cdots+n_k: n_1,\ldots, n_k).$ 

We find out later that the aforementioned immediate consequence 
is exactly the same as the formula of Yan and Wang \cite{YW12}. 
Nonetheless, the notion of $p(q: n_1,\ldots, n_k)$ 
can be used for each $q$ 
in which $K_{n_1,\ldots,n_k}$ has an equitable $q$-coloring 
and the proof presented here is much shorter.

\section{Main Result} 
We introduce the notion of $p(q: n_1,\ldots, n_k)$ which can be computed 
in  linear-time.   

\begin{Def}\label{d1}
Assume that $K_{n_1,\ldots,n_k}$ has an  equitable $q$-coloring,  
and $d$ is the minimum value not less than $\lceil (n_1+\cdots +n_k)/q \rceil$  
such that (i) there are distinct $i$ and $j$ in which $n_i$ and $n_j$ 
are not divisible by $d,$ or 
(ii) there is $n_j$ with  $n_j/\lfloor n_j/d\rfloor  > d+1.$ 
Define $p(q: n_1,\ldots, n_k)= \lceil n_1/d \rceil+\cdots+\lceil n_k/d \rceil.$  
\end{Def}

\begin{lem}\label{L1} 
Assume that $G=K_{n_1,\ldots,n_k}$ has an  equitable $q$-coloring.   
Then $G$ has an equitable $r$-coloring for each $r$ satisfying 
$p(q: n_1,\ldots, n_k) \leq r \leq q.$ 
\end{lem} 

\begin{Pf} 
Let $p=p(q: n_1, \ldots, n_k)$ and $N=n_1+\cdots +n_k.$ 
We prove by reverse induction that $G$ has an equitable $r$-coloring  
for each $r$ satisfying $p \leq r \leq q.$ 
By assumption, $G$ has an equitable  $q$-coloring. 
Consider $r$ in which $p < r \leq q$  and $G$ has an equitable $r$-coloring $f.$  
We show that $G$ has an equitable $(r-1)$-coloring. 
Let $b= \lceil N/r \rceil.$ 
By assumption, there are  integers $r_i$ and $s_i$ such that  
$f$ partitions $X_i$ into $r_i-s_i$ color classes of size $b$ and 
$s_i$ color classes of size $b-1$ where $r= r_1+\cdots+r_k.$ 
Thus $n_i = (r_i-s_i)b+s_i(b-1) = r_i b -s_i $ for each $i.$  

CASE 1: Some $j$ has $r_j \neq \lceil n_j/b \rceil.$ 
Note that $n_j = \lceil n_j/b \rceil b - g_j$ for some $g_j$  
satisfying $0 \leq g_j \leq b-1.$ 
Now, we have $r_j b -s_j = \lceil n_j/b \rceil b - g_j.$  
Thus $(r_j -\lceil n_j/b \rceil) b =s_j-g_j.$   
Combining with the fact $r_j \neq \lceil n_j/b \rceil, 0 \leq g_j \leq b-1,$ 
and  $s_j$ is positive,  
we have $s_j -g_j$ is a positive multiple of $b.$ 
From $n_j = (r_j-s_j)b+s_j(b-1),$ 
we can rewrite $n_j=(r_j-s_j+b-1)b+(s_j-b)(b-1).$ 
Since $s_j -r$ is a positive multiple of $b,$  we have $s_j-b$ is nonnegative. 
Thus we can partition $X_j$ into  $r_j-s_j+b-1$ color classes of size $b$ and 
$s_j-b$ color classes of size $b-1.$ 
That is, we can partition $X_j$ into $r_j-1$ color classes of size $b$ or $b-1.$ 
Since we can partition other $X_i$s into $r_i$ color classes of size $b$ or $b-1$ 
and 
$(\sum_{i \neq j}r_i)+(r_j-1)$ $=(\sum_{i=1}^k r_i)-1=r-1,$  
the graph $G$ has an equitable $(r-1)$-coloring. 

CASE 2: $r_i = \lceil n_i/b \rceil$ for each $i.$ 
Since $r > p$ and the condition of $d,$ we have $d > b.$ 
Thus $b$ violates conditions (i) and (ii) of $d$ in Definition \ref{d1}.  
Consequently, there are at least $k-1$ of $n_i$s which are a multiple of $b$ 
and  $n_j/\lfloor n_j/b\rfloor  \leq b+1$ for each $j.$ 
Without loss of generality, we assume $n_i=r_i b$ for each $i \geq 2.$ 

SUBCASE 2.1: $n_1 \neq r_1 b.$ 
Then  $b <  n_1/\lfloor n_1/b\rfloor =n_1/(\lceil n_1/b\rceil -1) = n_1/(r_1-1).$  
Since $b$ violates condition (ii), 
we have $n_1/(r_1-1)=n_1/\lfloor n_1/b\rfloor \leq b+1.$ 
Thus  $b<n_1/(r_1-1) \leq b+1.$ 
Consequently, we can partition $n_1$ into $r_1-1$ color classes of size $b$ or $b+1.$ 
Combining with $r_i$ color classes of $X_i$ of size $b$ for $i \geq 2,$ 
we have an equitable $(r-1)$-coloring.  

SUBCASE 2.2 $n_i = r_i b$ for each $i.$ If there is $j$ such that 
$n_j/(r_j-1) \leq b+1,$ then we have an equitable $(r-1)$-coloring as in subcase 2.1. 
Thus we assume further that $n_i/(r_i-1) > b+1$ for each $i.$ 
We claim that $b+1=d$ and $\lceil n_i/b\rceil = \lceil n_i/(b+1)\rceil = \lceil n_i/d\rceil.$ 
If the claim holds, we have $r =  \sum_{i=1}^k \lceil n/b \rceil=$ 
$\sum_{i=1}^k \lceil n/d \rceil = p$ which contradicts to the fact that $r>p.$ 
Thus this situation is impossible.  

To prove the claim, suppose to the contrary that $n_i$ is divisible by $b+1$ for some $i.$ 
Since $n_i = r_ib,$ we have $r_i =t_i(b+1)$ for some positive integer $t_i.$ 
Thus $n_i/(r_i-1)= t_i(b+1)b/(t_i(b+1)-1) = b+ b/(t_i(b+1)-1) \leq b+1$ 
which contradicts to the fact that $n_i/(r_i-1) > b+1.$ 
Thus $n_i$ is not divisible by $b+1$ for each $i.$  
Consequently, $b+1=d$ by condition (i). 
Since $n_i=r_ib$ and $n_i/(r_i-1) > b+1$ for each $i,$ 
we have $r_i = n_i/b > n_i/(b+1)> r_i-1.$ 
This leads to $r_i=\lceil n_i/b\rceil = \lceil n_i/(b+1)\rceil.$ 
Thus, we have the claim and this completes the proof.       
\end{Pf}

\begin{lem}\label{L2} 
Assume that $G=K_{n_1,\ldots,n_k}$ has an  equitable $q$-coloring    
and $p=p(q: n_1, \ldots, n_k).$
Then $G$ has no equitable $(p-1)$-coloring.
\end{lem} 
\begin{Pf}
Suppose to the contrary that $G$ has an equitable $(p-1)$-coloring.  
Then a partite set, say $X_1$ of size $n_1,$ is partitioned 
into at most $\lceil n_1/d \rceil -1$ color classes 
and a partite set $X_j$ of size $n_j$ is partitioned 
into at least $\lceil n_j/d \rceil$ color classes.  
Now we have at least color class containing vertices in $X_1$ 
with size at least $d+1.$ 
By (i) and (ii) in Definition \ref{d1}, we investigate 2 cases. 

CASE 1: there is some $X_j$ partitioned 
into at least $\lceil n_j/d \rceil+1$ color classes 
or there is some $n_j$ with $j\geq 2$ which is not divisible by $d.$ 
But then we have at least one color class containing vertices in $X_j$ 
with size at most $d-1.$ 
This contradicts to the fact the sizes of two color classes 
differ at most one. 

CASE 2: each $X_j$ with $j \geq 2$ 
has exactly $\lceil n_j/d \rceil$ color classes 
and $n_j$ is divisible by $d.$ 
Then $\lceil n_j/d \rceil=d$ for $j \geq 2.$  
Thus  $n_1$ has 
$n_1/\lfloor n_1/d\rfloor  > d+1$ by the condition (ii) of $d$ 
in  Definition \ref{d1}. 
But $X_1$ is partitioned into at least $\lceil n_1/d \rceil -1=\lfloor n_1/d\rfloor$ 
color classes.  
Thus we have at least one color class containing vertices in $X_1$ 
with size at least $d+2.$  
But each color class containing vertices in $X_j$ where $j\geq 2$ has 
size $d.$ Thus $G$ has no equitable $(p-1)$-coloring. 
\end{Pf}

From Lemmas \ref{L1} and \ref{L2}, we have the following theorem. 

\begin{Th}\label{T1}
Assume that $G=K_{n_1,\ldots,n_k}$ has an  equitable $q$-coloring.   
Then $p(q: n_1,\ldots, n_k)$ is the minimum $p$ such that 
$G$ has an equitable $r$-coloring for each $r$ satisfying $p \leq r \leq q.$ 
\end{Th}

\end{document}